\documentclass[submission%
]{dmtcs}

\usepackage{amssymb,amsmath,xspace,tocvsec2,mathrsfs,enumerate}
\usepackage[latin1]{inputenc}
\usepackage{soul,color}

\newtheorem{theorem}[equation]{Theorem}
\newtheorem{cor}[equation]{Corollary}
\newtheorem{lemma}[equation]{Lemma}
\newtheorem{proposition}[equation]{Proposition}

\newtheorem{remark}[equation]{Remark}

\numberwithin{equation}{section}

\newcommand{\commen}[1]{}
\newcommand{\one}[1]{\mathbf{1}_{#1 \times #1}}
\newcommand{\F}{\ensuremath{\mathbb{F}}}
\newcommand{\C}{\ensuremath{\mathbb{C}}}
\newcommand{\R}{\ensuremath{\mathbb{R}}}
\newcommand{\bp}{\ensuremath{\mathcal{P}}}
\newcommand{\calp}{\ensuremath{\mathcal{P}}}
\newcommand{\cals}{\ensuremath{\mathcal{S}}}
\newcommand{\ba}{\mathbf{a}}
\newcommand{\bc}{\mathbf{c}}
\newcommand{\bd}{\mathbf{d}}
\newcommand{\bn}{\mathbf{n}}
\newcommand{\bu}{\mathbf{u}}
\newcommand{\bv}{\mathbf{v}}
\newcommand{\diag}{\mathop{\mathrm{diag}}\nolimits}
\newcommand{\Gr}{\mathop{\mathrm{Gr}}\nolimits}

\author{Alexander Belton\addressmark{1}\thanks{Email: \email{a.belton@lancaster.ac.uk}} ,
  Dominique Guillot\addressmark{2}\thanks{Email: \email{dguillot@udel.edu}} ,
  Apoorva Khare\addressmark{3}\thanks{Email: \email{khare@stanford.edu}} \\ \and
  Mihai Putinar\addressmark{4,5}\thanks{Email: \email{mputinar@math.ucsb.edu, mihai.putinar@ncl.ac.uk}}}
\address{
	\addressmark{1}Lancaster University, Lancaster, UK \\
	\addressmark{2}University of Delaware, Newark, DE, USA \\
	\addressmark{3}Stanford University, Stanford, CA, USA \\
	\addressmark{4}University of California at Santa Barbara, CA, USA \\
	\addressmark{5}Newcastle University, Newcastle upon Tyne, UK
}

\title[Schur polynomials and matrix positivity
preservers]{Schur polynomials and\\ matrix positivity preservers}
\begin{document}
\maketitle

\begin{abstract}
\paragraph{Abstract.}
A classical result by Schoenberg (1942) identifies all real-valued
functions that preserve positive semidefiniteness (psd) when applied
entrywise to matrices of arbitrary dimension. Schoenberg's work has
continued to attract significant interest, including renewed recent
attention due to applications in high-dimensional statistics.
However, despite a great deal of effort in the area, an effective
characterization of entrywise functions preserving positivity in a fixed
dimension remains elusive to date. As a first step, we characterize new
classes of polynomials preserving positivity in fixed dimension. 
The proof of our main result is representation theoretic, and employs
Schur polynomials. An alternate, variational approach also leads to
several interesting consequences including
(a) a hitherto unexplored Schubert cell-type stratification of the cone
of psd matrices,
(b) new connections between generalized Rayleigh quotients of Hadamard
powers and Schur polynomials,
and
(c) a description of the joint kernels of Hadamard powers.
\medskip

\paragraph{R\'esum\'e.}
Un r\'esultat classique de Schoenberg (1942) fournit une
caract\'erisation des fonctions r\'eelles pr\'eservant la positivit\'e
lorsque appliqu\'ees aux entr\'ees des matrices semid\'efinie positives
de dimension arbitraire. Le travail de Schoenberg est toujours
d'actualit\'e, et a r\'ecemment re\c cu beaucoup d'attention suite \`a
ses applications aux statistiques de haute dimension. N\'eanmoins,
l'obtention d'une caract\'erisation utile des fonctions pr\'eservant la
positivit\'e lorsque la dimension est fixe demeure un probl\`eme ouvert.
Afin d'attaquer ce probl\`eme, nous caract\'erisons de nouvelles classes
de polyn\^omes qui pr\'eservent la positivit\'e en dimension finie. La
preuve de notre r\'esultat principal implique plusieurs id\'ees provenant
de la th\'eorie de la repr\'esentation, et utilise les polyn\^omes de
Schur. Nous explorons aussi une approche variationnelle parall\`ele qui
m\`ene \`a de nombreux r\'esultats int\'eressants:
(a) une stratification du c\^one des matrices semid\'efinies positives,
(b) de nouvelles connexions entre les quotients de Rayleigh
g\'en\'eralis\'es des puissances d'Hadamard et les polyn\^omes de Schur,
et (c) une description du noyau simultan\'e des puissances d'Hadamard.

\end{abstract}

\section{Introduction and main result}

Endomorphisms of matrix spaces with positivity constraints
have long been studied in connection with a variety of topics: the
geometry of classical domains in complex space, matrix monotone functions
\cite{Loewner34}, positive definite functions \cite{BCR-semigroups,
Berg-Porcu, Bochner-pd, Schoenberg42, vonNeumann-Schoenberg}, hyperbolic
or positive definite polynomials and global optimization algorithms
\cite{SIAM-optimization, HMPV}.
In this paper, we study the entrywise calculus on the cone of positive
semidefinite matrices, with the aim of characterizing positivity
preservers in that setting.

Given $\rho \in (0, \infty)$, let $D(0, \rho)$ and $\overline{D}(0,
\rho)$ denote the open and closed complex discs of radius $\rho$ centered
at the origin, respectively.
Given integers $1 \leq k \leq N$ and a set $I \subset \C$, let
$\bp_N^k(I)$ denote the set of positive semidefinite $N \times N$
matrices, with entries in $I$ and rank at most $k$. Let $\bp_N(I) :=
\bp_N^N(I)$. A function $f : I \to \C$ induces an entrywise map of matrix
spaces, sending $A = (a_{jk}) \in \bp_N(I)$ to $f[A] := (f(a_{jk}))$.
Starting from positive definite functions \cite{Bochner-pd, Schoenberg42,
vonNeumann-Schoenberg}, it is natural to classify all entrywise functions
$f[-]$ preserving positive semidefiniteness (positivity). It is an easy
consequence of the Schur product theorem \cite{Schur1911} that if $f :
(-\rho,\rho) \to \R$ is analytic with non-negative Taylor coefficients,
then $f[A] \in \bp_N$ for all $A \in \bp_N$ and all $N \geq 1$. A
celebrated result of Schoenberg shows the converse.

\begin{theorem}[Schoenberg, \cite{Schoenberg42}]\label{Tschoenberg}
Given a continuous function $f : [ {-1}, 1 ] \to \R$, the entrywise map
$f[ - ] : \bp_N( [ {-1}, 1 ] ) \to \bp_N( \R )$ for all $N \geq 1$ if and
only if $f$ is analytic on $[ {-1}, 1 ]$ and absolutely monotonic on
$[0,1]$, i.e., $f$ has a Taylor series with non-negative coefficients
convergent on $\overline{D}( 0, 1 )$.
\end{theorem}

Schoenberg's theorem and its ramifications were persistently examined and
revisited, see for instance Rudin~\cite{Rudin59},
Berg--Christensen--Ressel--Porcu~\cite{BCR-semigroups, Berg-Porcu,
Christensen_et_al78}, Hiai~\cite{Hiai2009}, to cite only a few. The
present investigation evolves out of Schoenberg's result by imposing the
challenging condition of dealing with matrices of \textit{fixed}
dimension. This is a much harder question, that is open despite
tremendous activity in the field.

It is worth recalling that Schoenberg was motivated by the problem of
isometrically embedding positive definite metrics into Hilbert space; see
e.g.~\cite{vonNeumann-Schoenberg}. In~\cite{Schoenberg42}, he sought to
classify positive definite functions on spheres $S^{d-1} \subset \R^d$.
This can be reformulated via Gram matrices, as classifying the entrywise
functions preserving positivity on correlation matrices of all
dimensions, with rank at most $d$.
A strong need to study the fixed dimension case also arises out of
current demands from the fast expanding field of data science.
In modern settings, functions $f$ are often applied entrywise to
high-dimensional correlation matrices $A$, in order to improve their
properties (better conditioning, Markov random field structure, etc.);
see e.g.~\cite{bickel_levina, hero_rajaratnam, Rothman_et_al_JASA}. 
The ``regularized'' matrices $f[A]$ are ingredients in further
statistical procedures, for which it is critical that they be positive
semidefinite. Also, in applications the dimension of the
problem is known, and so, preserving positivity in all
dimensions unnecessarily limits the class of functions that can be used.
There is thus strong motivation from applications to study the fixed
dimension case.

While characterization results have recently been obtained in fixed
dimension under additional rank and sparsity constraints arising in
practice \cite{GKR-lowrank, Guillot_Khare_Rajaratnam2012,
Guillot_Rajaratnam2012b}, the original problem in fixed dimension has
remained open for more than 70 years. A necessary condition for
continuous functions was developed by Horn (and attributed to Loewner) in
his doctoral thesis \cite{horn}. The result was recently extended in
\cite{GKR-lowrank} to low-rank matrices, and without the continuity
assumption:

\begin{theorem}[{Horn \cite{horn}, Guillot--Khare--Rajaratnam
\cite{GKR-lowrank}}]\label{Thorn}
Suppose $f : I \to \R$, where $I := ( 0, \rho )$ and
$0 < \rho \leq \infty$. Fix an integer $N \geq 2$ and suppose that
$f[ A ] \in \bp_N( \R )$ for any $A \in \bp_N^2( I )$ of the form
$A = a \one{N} + \bu \bu^T$, where $a \in ( 0, \rho )$,
$\bu \in [ 0, \sqrt{\rho - a} )^N$, and $\one{N} \in \bp_N^1(\R)$ has
entries all one. Then $f \in C^{N - 3}(I)$, with
$f^{( k )}( x ) \geq 0$ for all $x \in I, \ 0 \leq k \leq N - 3$,
and $f^{(N - 3)}$ is a convex non-decreasing function on~$I$. If,
further, $f \in C^{N - 1}( I)$, then $f^{( k )}( x ) \geq 0$ for all
$x \in I$ and $0 \leq k \leq N - 1$.
\end{theorem}

Note that all real power functions $x^\alpha$ preserve positivity on
$\bp_N^1( ( 0, \rho ) )$, yet such functions need not have even a single
positive derivative on $( 0, \rho )$. However, Theorem~\ref{Thorn} shows
that working with a small one-parameter extension of $\bp_N^1( ( 0, \rho
) )$ guarantees that $f^{( k )}$ is non-negative on $( 0, \rho )$ for $0
\leq k \leq N - 3$.

Theorem~\ref{Thorn} is sharp, since the entrywise power $x^\alpha$, for
$\alpha \in (N-2, N-1)$, preserves positivity on $\bp_N( ( 0, \rho ) )$,
but not on $\bp_{N + 1}( ( 0, \rho ) )$.
See~\cite{FitzHorn, GKR-crit-2sided, Hiai2009} for more on entrywise
powers preserving positivity. Consequently, in this paper we study
analytic functions which preserve $\bp_N$ for fixed $N$, when applied
entrywise.
Note that any analytic function mapping $( 0, \rho )$ to $\R$ necessarily
has real Taylor coefficients. Now a variant of Theorem~\ref{Thorn} for
analytic functions, obtained using generalized Vandermonde matrices, shows
that the same conclusions hold if one works merely with rank-one
matrices:

\begin{lemma}\label{Lpos-coeff}
Let $0 < \rho \leq \infty$ and $f( z ) = \sum_{k \geq 0} c_k z^k : D( 0,
\rho ) \to \R$ be analytic. If $f[ - ] : \bp_N^1( ( 0, \rho )) \to
\bp_N(\R)$ for some integer $N \geq 1$, then the first $N$ non-zero
Taylor coefficients $c_j$ are strictly positive.
\end{lemma}

Given $f( z ) = \sum_{k \geq 0} c_k z^k$ such that $c_0$, \ldots,
$c_{N - 1} > 0$, a natural challenging question to ask is if the next
non-zero coefficient $c_M$ can be negative; and if so, to provide a
negative threshold for the coefficient $c_M$, where $M \geq N$.
Resolving these questions, open since Horn's 1969 paper, provides a
quantitative version of Schoenberg's theorem. Our main result answers
these questions in the affirmative, and illustrates the complexity of
the negative threshold bound. It is also surprising that preserving
positivity on $\bp_N(\overline{D}( 0, \rho ) )$ is equivalent to
preserving positivity on the much smaller set of real rank-one
matrices, $\bp_N^1( ( 0, \rho ) )$.

\begin{theorem}\label{Tthreshold}
Fix $\rho > 0$ and integers $N \geq 1$, $M \geq 0$ and let
$f( z ) := \sum_{j = 0}^{N - 1} c_j z^j + c' z^M$ be a
polynomial with real coefficients. Then the following are equivalent.
\begin{enumerate}
\item $f[ - ]$ preserves positivity on
$\bp_N( \overline{D}( 0, \rho ) )$.

\item Either
$c_0$, \ldots, $c_{N - 1}$, $c' \geq 0$, or
$c_0$, \ldots, $c_{N - 1} > 0$ and
$c' \geq -\mathcal{C}( \bc; z^M; N, \rho )^{-1}$, where
\begin{equation}
\bc := ( c_0, \ldots, c_{N - 1} ), \qquad
\mathcal{C}( \bc; z^M; N, \rho ) := %
\sum_{j = 0}^{N - 1} \binom{M}{j}^2 \binom{M - j - 1}{N - j - 1}^2 %
\frac{\rho^{M - j}}{c_j}.
\end{equation}

\item $f[ - ]$ preserves positivity on $\bp_N^1( ( 0, \rho ) )$.
\end{enumerate}
\end{theorem}

Notice that the condition $c_0$, \ldots, $c_{N - 1} \geq 0$ follows
from Lemma~\ref{Lpos-coeff}. Theorem \ref{Tthreshold} now
provides the first construction of a polynomial that preserves positivity
on $\bp_N$, but not on $\bp_{N+1}$. Indeed, this is the case when
$-\mathcal{C}( \bc; z^M; N, \rho)^{-1} \leq c_N < 0$, by Theorem
\ref{Thorn}.

\begin{remark}
Theorem~\ref{Tthreshold} can naturally be used to provide a sufficient
condition for an arbitrary analytic function to preserve positivity on
$\bp_N( \overline{D}( 0, \rho ) )$. The reader is referred to
\cite{BGKP-fixeddim} for more details.
\end{remark}

\section{Proof of the main result}

We now sketch the proof of Theorem~\ref{Tthreshold}.
Recall that the Schur product theorem provides the first examples of
entrywise functions preserving positivity, namely, the functions of the
form $\sum_{k=0}^\infty c_k z^k$ with $c_k \geq 0$. That these are the
only functions preserving positivity in all dimensions is Schoenberg's
theorem (Theorem \ref{Tschoenberg}).
In some sense, our proof of the fixed dimension case in
Theorem~\ref{Tthreshold} returns to Schur by crucially using symmetric
functions among other techniques, specifically, Schur polynomials and
Schur complements. Indeed, the technical heart of the proof is an
explicit Jacobi--Trudi type identity, which is valid in any field and
may be interesting in its own right.

Given a partition, i.e., a non-increasing $N$-tuple of
non-negative integers $\bn = ( n_N \geq \cdots \geq n_1 )$, define the
corresponding \emph{Schur polynomial} $s_\bn( x_1, \ldots, x_N )$ over a
field $\F$ with at least $N$ elements, to be the unique polynomial
extension to $\F^N$ of
$\displaystyle s_\bn( x_1, \ldots, x_N ) := %
\frac{\det ( x_i^{n_j + N - j} )}{\det ( x_i^{N - j} )}$
for pairwise distinct $x_i \in \F$. Note that the denominator
equals the Vandermonde determinant
$\Delta_N( x_1, \ldots, x_N ) := \det ( x_i^{N - j} ) = %
\prod_{1 \leq i < j \leq N} ( x_i - x_j )$; thus,
\begin{equation}\label{schureval}
s_\bn( 1, \ldots, z^{N - 1} ) = %
\prod_{1 \leq i < j \leq N} %
\frac{z^{n_j + j} - z^{n_i + i}}{z^j - z^i}, \qquad %
s_\bn( 1, \ldots, 1 ) = \prod_{1 \leq i < j \leq N} %
\frac{n_j - n_i + j - i}{j - i}.
\end{equation}
The last equation can also be deduced from the Weyl Character Formula in
type~$A$; see, for example, \cite[Chapter~I.3, Example~1]{Macdonald} for
more details about Schur polynomials and the theory of symmetric
functions.

To prove Theorem~\ref{Tthreshold}, we study the determinants of a linear
pencil
\[
p(t) = p_t[ A ] := \det \left( t ( c_0 {\bf 1}_{N \times N} + c_1 A +
\cdots + c_{N - 1}A^{\circ ( N - 1 )} ) - A^{\circ M} \right)
\]
for a general rank-one matrix $A = \bu \bv^T$, where
$\bu = ( u_1, \ldots, u_N )^T$ and $\bv := ( v_1, \ldots, v_N )^T
\in \F^N$ for $N \geq 1$. The technical heart of the proof involves
the following explicit determinantal identity.

\begin{theorem}\label{Pjacobi-trudi}
Let $M \geq N \geq 1$ be integers, and $c_0$, \ldots, $c_{N - 1} \in
\F^\times$ be non-zero scalars in any field $\F$.
Define the polynomial
$p_t( z ) := t ( c_0 + \cdots + c_{N - 1} z^{N - 1} ) - z^{M}$.
Now define the hook partition
\begin{equation}
\mu( M, N, j ) := ( M - N + 1, 1, \ldots, 1, 0, \ldots, 0 )
\end{equation}
($N-j-1$ ones, $j$ zeros) for $0 \leq j < N$.
Then the following identity holds for all ${\bf u}, {\bf v} \in
\F^N$:
\begin{equation}\label{Ejacobi-trudi}
\det p_t[ {\bf u} {\bf v}^T ] = t^{N - 1} \Delta_N( {\bf u} ) %
\Delta_N( {\bf v} ) \prod_{j = 0}^{N - 1} c_j %
\Bigl( t - %
\sum_{j = 0}^{N - 1} \frac{s_{\mu( M, N, j )}( {\bf u} ) %
s_{\mu( M, N, j )}( {\bf v} )}{c_j} \Bigr).
\end{equation}
Moreover,
$\displaystyle s_{\mu( M, N, j )}( 1, \ldots, 1 ) = %
\binom{M}{j} \binom{M-j-1}{N-j-1}$ for all $0 \leq j < N$.
\end{theorem}

\noindent {\bf Sketch of proof.}
We first show the following fact:
\textit{Let $A := \bu \bv^T$ for $\bu$, $\bv \in \F^N$.
Given a strict partition
$\bn = ( n_m > n_{m-1} > \cdots > n_1 )$
and scalars $( c_{n_1}, \ldots, c_{n_m} ) \in \F^m$,
the following determinantal identity holds:
\begin{equation}\label{Ecauchy-binet1}
\det \sum_{j = 1}^m c_{n_j} A^{\circ n_j} = %
\Delta_N( \bu ) \Delta_N( \bv ) \sum_{\bn' \subset \bn, \ | \bn'| = N}
s_{\lambda( \bn' )}( \bu ) s_{\lambda( \bn' )}( \bv ) %
\prod_{k = 1}^N c_{n'_k}.
\end{equation}
Here,
$\lambda( \bn') := %
( n'_N - N + 1 \geq n'_{N-1} - N + 2 \geq \cdots \geq n'_1 )$
is obtained by subtracting the staircase partition
$( N - 1, \ldots, 0 )$ from $\bn' := ( n'_N > \cdots > n'_1 )$,
and the sum is over all subsets $\bn'$ of cardinality~$N$.}\medskip

The proof of~\eqref{Ecauchy-binet1} uses the matrix
$X( \bu, \bn ) := ( u_j^{n_k} )_{1 \leq j \leq N, 1 \leq k \leq m}$
and the Cauchy--Binet formula applied to
$\displaystyle \sum_{j = 1}^m c_{n_j} A^{\circ n_j} = %
X( \bu, \bn ) \cdot \diag( c_{n_1}, \ldots, c_{n_m} ) \cdot %
X( \bv, \bn )^T$.\medskip

Using~\eqref{Ecauchy-binet1}, we now prove~\eqref{Ejacobi-trudi}.
Recall the Laplace formula:
if $B$ and $C$ are $N \times N$ matrices, then
$\displaystyle \det( B + C ) = %
\sum_{\bn \subset \{ 1, \ldots, N \}} \det M_\bn( B; C)$,
where $M_\bn( B; C )$ is the matrix formed by replacing the rows of
$B$ labelled by elements of $\bn$ with the corresponding rows of $C$.
In particular, if $B = \sum_{j = 0}^{N - 1} c_j A^{\circ j}$ then
\begin{equation}
\det p_t[ A ] = \det( t B - A^{\circ M} ) = t^N \det B - %
t^{N - 1} \sum_{j = 1}^N \det M_{\{ j \}}( B; A^{\circ M} ),
\end{equation}

\noindent since the determinant in each of the remaining terms contains
at least two rows of the rank-one matrix $A^{\circ M}$. 
Applying~\eqref{Ecauchy-binet1} yields:
$\det B = \Delta_N( \bu ) \Delta_N( \bv ) c_0 \cdots c_{N - 1}$.
Moreover, the coefficient of $t^{N-1}$ is precisely $\det p_1[A] - \det
B$, and $\det p_1[A]$ can be computed using~\eqref{Ecauchy-binet1} again,
to yield:
\[
\det p_1[ A ] = \det B - \Delta_N( \bu ) \Delta_N( \bv ) %
c_0 \cdots c_{N - 1} \sum_{j = 0}^{N - 1} %
\frac{s_{\mu( M, N, j )}( \bu ) s_{\mu( M, N, j)}( \bv )}{c_j},
\]
since
$\mu( M, N, j ) = %
\lambda\bigl( ( M, N - 1, N - 2, \ldots, j + 1, %
\widehat{j}, j - 1, \ldots, 0 ) \bigr)$ for $0 \leq j < N$.
This proves the identity~\eqref{Ejacobi-trudi}.
The final assertion follows from~\eqref{schureval}, or by
using the dual Jacobi--Trudi (Von N{\"a}gelsbach--Kostka) identity
\cite[Chapter I, Eqn.~(3.5)]{Macdonald}, or Stanley's hook-content formula
\cite[Theorem 15.3]{Stanley}.\qed\medskip

Equipped with Theorem~\ref{Pjacobi-trudi}, we now outline how to show the
main result.\medskip

\noindent {\bf Sketch of proof of Theorem~\ref{Tthreshold}.}
Clearly $(1) \implies (3)$.

Now if $0 \leq M < N$,
then $\mathcal{C}( \bc; z^M; N, \rho ) = c_M^{-1}$. Thus, in this
case, $(2) \implies (1)$ by the Schur product theorem, and moreover,
$(3) \implies (2)$ by Lemma~\ref{Lpos-coeff}.

Suppose for the remainder of the proof that $M \geq N$; we also set $c_M
:= c'$. We first show that $(3) \implies (2)$. It suffices to consider
the case $c_M < 0 < c_0$, \ldots, $c_{N - 1}$. Define $p_t( z )$ as in
Theorem~\ref{Pjacobi-trudi}, and set $t := |c_M|^{-1}$. By
Equation~\eqref{Ejacobi-trudi},
\begin{equation}\label{Enecessary}
0 \leq \det p_t[ \bu \bu^T ] = t^{N - 1} \Delta_N( \bu )^2 %
c_0 \cdots c_{N - 1} \Bigl( t - \sum_{j = 0}^{N - 1} %
\frac{s_{\mu( M, N, j )}( \bu )^2}{c_j} \Bigr).
\end{equation}
Set $u_k := \sqrt{\rho} ( 1 - t' \epsilon_k )$, with pairwise distinct
$\epsilon_k \in ( 0, 1 )$ , and $t' \in ( 0, 1 )$. Thus,
$\Delta_N( \bu ) \neq 0$. Taking the limit as $t' \to 0^+$, since
the final term in~\eqref{Enecessary} must be non-negative, it follows by
Theorem~\ref{Pjacobi-trudi} that
\[ t = |c_M|^{-1} \geq \sum_{j = 0}^{N - 1} %
\frac{s_{\mu( M, N, j )}( \sqrt{\rho}, \ldots, \sqrt{\rho} )^2}{c_j} %
 = \sum_{j = 0}^{N - 1} s_{\mu( M, N, j )}( 1, \ldots, 1 )^2 %
\frac{\rho^{M - j}}{c_j}
 = \mathcal{C}( \bc; z^M; N, \rho ). \]

It remains to show that $(2) \implies (1)$ when $M \geq N$ and
$c_M < 0 < c_0$, \ldots, $c_{N - 1}$. The major step involves showing
that $f[-]$ preserves positivity on
$\bp_N^1( \overline{D}( 0, \rho ) )$. Given $1 \leq m \leq N$, define
\begin{equation}
C_m := \sum_{j = 0}^{m - 1} s_{\mu( M-N+m, m, j )}( 1, \ldots, 1 )^2 %
\frac{\rho^{m + M-N - j}}{c_{N - m + j}}
= \mathcal{C}( \bc_m; z^{M - N + m}; m, \rho ),
\end{equation}
where $\bc_m := ( c_{N - m}, \ldots, c_{N - 1} )$. One now shows that
\[
0 < C_1 = \rho^{M-N+1}/c_{N-1} < C_2 < \cdots < C_N = \mathcal{C}( \bc;
z^M; N, \rho ).
\]

Next, we claim that for all $1 \leq m \leq N$ and $A = \bu \bu^*
\in \bp_N^1( \overline{D}( 0, \rho ) )$, every principal $m \times m$
submatrix of the matrix
\begin{equation}
L := C_m ( c_{N - m} \one{N} + c_{N - m + 1} A + \cdots + %
c_{N - 1} A^{\circ ( m - 1 )} ) - A^{\circ ( m + M - N )}
\end{equation}
is positive semidefinite. Notice that the rank-one case of~(1) follows by
setting $m = N$.

The claim is shown by induction on $m$, with the $m=1$ case immediate.
Suppose the result holds for $m-1 \geq 1$. Henceforth, given a non-empty set
$\bn \subset \{ 1, \ldots, N \}$ and an $N \times N$ matrix $B$,
denote by $B_\bn$ the principal submatrix of $B$ consisting of those rows
and columns labelled by elements of $J$; similarly define the subvector
$\bu_\bn$ of a vector $\bu$. If $\bn \subset \{ 1, \ldots, N \}$ has
cardinality $m$ then,
\[
\det L_\bn = %
C_m^{m - 1} | \Delta_m( \bu_\bn )|^2 \prod_{j = 1}^m c_{N - j} %
\Bigl( C_m - \sum_{j = 0}^{m - 1} %
\frac{| s_{\mu( M - N + m, m, j )}( \bu_\bn ) |^2}{c_{N - m + j}} \Bigr)
\]

\noindent by Theorem~\ref{Pjacobi-trudi} with $\bv = \overline{\bu_n}$.
By the triangle inequality in $\C$ and the fact that the
coefficients of Schur polynomials are non-negative, one shows that
$\det L_\bn \geq 0$ if $|\bn| = m$. If on the other hand $|\bn| = k<m$,
\begin{align*}
L_\bn & \geq C_m ( c_{N - k} A_\bn^{\circ ( m - k )} + \cdots + %
c_{N - 1} A_\bn^{\circ ( m - 1 )} ) - A_\bn^{\circ ( m + M - N)} \\
 & \geq A_\bn^{\circ ( m - k )} \circ \bigl( C_k %
( c_{N - k} \one{k} + \cdots + c_{N - 1} A_\bn^{\circ ( k - 1 )} ) %
- A_\bn^{\circ ( k + M - N )} \bigr),
\end{align*}
since $C_m > C_k$. Hence by the induction hypothesis, all principal $m
\times m$ submatrices of $L$ are positive semidefinite. This shows the
claim by induction on $m$, whence $f[-]$ preserves positivity on
$\bp_N^1( \overline{D}( 0, \rho ) )$.

Finally, to show the result for matrices in $\bp_N$ of all ranks, we
induct on $N \geq 1$. Suppose (1) holds for $N-1 \geq 1$, and define
$p_t[ B; M, \bd ] := t ( d_0 {\bf 1} + d_1 B + \cdots + %
d_{n - 1} B^{\circ ( n - 1 )}) - B^{\circ ( n + M )}$
for any square matrix $B$, where $t$ is a real scalar and $n$ is the
length of the tuple $\bd = ( d_0, \ldots, d_{n - 1} )$. It suffices to
show the claim that $p_t[ A; M - N, \bc ] \geq 0$ for
all~$t \geq \mathcal{C}( \bc; z^M; N, \rho )$ and all
$A = ( a_{i j} ) \in \bp_N( \overline{D}( 0, \rho ) )$.

To show the claim, define $\bu := (a_{iN} / \sqrt{a_{NN}})^T \in \C^N$,
and use \cite[Lemma~2.1]{FitzHorn} to show that
\begin{equation}\label{eqn:convex}
p_t[ A; M - N, \bc ] = p_t[ \bu \bu^*; M - N, \bc ]
+ \int_0^1 ( A - \bu \bu^* ) \circ %
M p_{t / M}[ \lambda A + ( 1 - \lambda ) \bu \bu^*; M - N,
\bc' ] d\lambda,
\end{equation}
where the $(N - 1)$-tuple
$\bc' := (c _1, 2 c_2, \ldots, ( N - 1 ) c_{N - 1} )$.
Notice that $A - \bu \bu^*$ is the padding by a zero row and column,
of the Schur complement of $a_{NN}$ in $A$. Therefore the integrand
in~\eqref{eqn:convex} is positive semidefinite if the matrix
$p_{t / M}[ A_\lambda; M - N, \bc' ]$ is, where
$A_\lambda \in \bp_{N - 1}( \overline{D}( 0, \rho ) )$ is obtained by
deleting the final row and column of
$\lambda A + ( 1 - \lambda ) \bu \bu^*$. Finally, that
$p_{t / M}[ A_\lambda; M - N, \bc' ]$ is positive semidefinite follows by
computing that
$\mathcal{C}( \bc; z^M; N, \rho ) \geq
M \mathcal{C}( \bc'; z^{M - 1}; N - 1, \rho )$.\qed\medskip

\section{Consequences of the main theorem}

Theorem~\ref{Tthreshold} leads to a host of consequences that initiate
the development of an entrywise matrix calculus, in parallel to the
well-studied functional calculus. We now discuss two of these
consequences in detail: linear matrix inequalities and connections to
Rayleigh quotients.

\subsection{Linear matrix inequalities for Hadamard powers}

Theorem~\ref{Tthreshold} can be equivalently reformulated as a linear
matrix inequality that controls the spectrum of linear combinations of
Hadamard powers of $A$.

\begin{theorem}
Fix $\rho > 0$, integers $M \geq N \geq 1$, and scalars $c_0$, \ldots,
$c_{N - 1} > 0$. Then
\begin{equation}\label{Elmi}
A^{\circ M} \leq \mathcal{C}( \bc; z^M; N, \rho ) \cdot \bigl( %
c_0 \one{N} + c_1 A + \cdots + c_{N-1} A^{\circ (N - 1)} \bigr),
\qquad \forall A \in \bp_N( \overline{D}( 0, \rho) ),
\end{equation}

\noindent where $\leq$ stands for the Loewner ordering. Moreover, the
constant $\mathcal{C}( \bc; z^M; N, \rho )$ is sharp in~\eqref{Elmi}.
\end{theorem}

\noindent Notice here that the right-hand side of~\eqref{Elmi} cannot
involve fewer Hadamard powers, by Lemma~\ref{Lpos-coeff}.

A refined analysis of the proof of Theorem~\ref{Tthreshold} shows that
the matrix $f[ A ]$ is generically positive definite, in a strong sense:

\begin{theorem}
Fix $\rho > 0$, integers $M \geq N \geq 1$, and scalars $c_0, \ldots,
c_{N - 1} > 0$.
\begin{enumerate}
\item Suppose $N>1$, and $A \in \bp_N( \overline{D}( 0, \rho ) )$ has a
row or column with pairwise distinct entries. Define $f(z) := c_0 +
\cdots + c_{N-1} z^{N-1} - \mathcal{C}( \bc; z^M; N, \rho )^{-1} z^M$.
Then $f[A]$ is positive definite.

\item In particular, equality in~\eqref{Elmi} is never attained on
$\bp_N( \overline{D}( 0, \rho) )$ unless $(N,A) = (1,\rho)$.
\end{enumerate}
\end{theorem}

\noindent As with the main theorem, the proof of (1) and (2) crucially
uses symmetric functions, specifically, connections between Schur
polynomials and Young tableaux.\medskip

\noindent {\bf Sketch of proof.}
It is also easy to show (2) for $N=1$.
Thus, assume $N>1$. We first show that~(2) follows from~(1).
Indeed, if~(1) holds, then using $\bu \in [0, \sqrt{\rho}]^N$ with
distinct entries, it follows that $f : [0,\rho] \to (0,\infty)$. Thus,
$f[A]$ has positive diagonal entries for any $A \in \bp_N( \overline{D}(
0, \rho  ) )$, whence~\eqref{Elmi} is never an equality.

It remains to prove~(1); here we sketch only the argument for $A = \bu
\bu^*$ of rank one. Thus $\bu$ has pairwise distinct entries by
assumption, with $N>1$. Suppose for contradiction that $\det f[A] = 0$.
Then,
\[
\sum_{j = 0}^{N - 1} \frac{| s_{\mu( M, N, j )}%
( \sqrt{\rho}, \ldots, \sqrt{\rho} )|^2}{c_j} =: %
\mathcal{C}( \bc; z^M; N, \rho ) = %
\sum_{j = 0}^{N - 1} \frac{| s_{\mu( M, N, j )}( \bu ) |^2}{c_j}
\]

\noindent using~\eqref{Ejacobi-trudi}. Now use the fact that
the coefficients of any Schur polynomial are non-negative, to show that
\[
| s_{\mu( M, N, j )}( \bu ) | = %
s_{\mu( M, N, j )}( \sqrt{\rho}, \ldots, \sqrt{\rho} )
= \binom{M}{j} \binom{M - j - 1}{N - j - 1} \rho^{( M - j ) / 2} \qquad
\forall j.
\]

Consider the case $j = N - 1$, which corresponds to the partition
$\mu( M, N, N - 1 ) = ( M - N + 1, 0, \ldots, 0 )$. By
\cite[Chapter~I, Equation~(5.12)]{Macdonald}, the Schur
polynomial $s_{\mu( M, N, N - 1 )}$ is a sum of $\binom{M}{N-1}$
monomials $\bu^{\bf t}$ corresponding to semi-standard Young tableaux
(i.e., ${\bf t} = (t_1, \ldots, t_N)$ with $\sum_{j=0}^{N-1} t_j =
M-N+1$). Now using the triangle inequality, all monomials $\bu^{\bf t}$
are equal; since $u_1^{M-N} u_j$ is such a monomial for all $j$, we get
$u_1 = \cdots = u_N$, a contradiction. Hence $\det f[A] > 0$ as claimed.
\qed\medskip

Theorem~\ref{Tthreshold} also fits naturally into the framework of
spectrahedra and the matrix cube problem~\cite{SIAM-optimization,
Nemirovski}; see~\cite{BGKP-fixeddim} for more details.

\subsection{Rayleigh quotients}\label{Srayleigh}

Given a domain $K \subset \C$,
functions $g, h : K \to \C$, and a set of matrices $\displaystyle \calp
\subset \cup_{N \geq 1} \bp_N(K)$, define $\mathcal{C}( h; g; \calp)$
to be the smallest real number such that
$g[A] \leq \mathcal{C}( h; g; \calp ) \cdot h[A]$ for all $A \in \calp$.
That is, $\mathcal{C}( h; g; \calp )$ is the \textit{extreme
critical value} of the family of linear pencils
$\{ -g[A] + \R h[A] : A \in \calp \}$. This notation helps achieve a
uniform and consistent formulation of the aforementioned theorems by
Schoenberg and Horn, Theorem~\ref{Tthreshold} and its consequences, as
well as other results in the literature. See \cite[Section
6]{BGKP-fixeddim} for a comprehensive survey of numerous such results.

Given $\bc = (c_0, \ldots, c_{N-1}) \in (0,\infty)^N$, define the
polynomial $h_\bc(z) := \sum_{j=0}^{N-1} c_j z^j$. By
Theorem~\ref{Tthreshold},
\begin{equation}
\mathcal{C}( \bc; z^M; N, \rho) = \mathcal{C}( h_\bc; z^M; \calp ),
\qquad \forall \bp_N^1((0,\rho)) \subset \calp \subset
\bp_N( \overline{D}( 0, \rho ) ).
\end{equation}

We now discuss an alternate, variational approach to proving
Theorem~\ref{Tthreshold}, which proceeds as follows:
\begin{enumerate}
\item[(I)] Bound $A^{\circ M}$ by lower Hadamard powers for a
\textit{single} matrix $A$, i.e., by $\alpha_A \cdot h_\bc[A]$ for the
smallest constant $\alpha_A > 0$.

\item[(II)] Now take the supremum of $\alpha_A$ over all matrices $A \in
\bp_N( \overline{D}( 0, \rho) )$.
\end{enumerate}

\noindent Notice that the first step (I) simply involves computing the
extreme critical value $\alpha_A = \mathcal{C}( h_\bc; z^M; A)$, using
the above notation. This and an improved understanding of $\ker
h_\bc[A]$, can be achieved as follows:

\begin{proposition}\label{Crayleigh}
Fix $\rho > 0$, integers $M \geq N \geq 1$, scalars $c_0$, \ldots, $c_{N
- 1} > 0$, and $A \in \bp_N( \C ) $. Define
\begin{equation}\label{Esimult}
\mathcal{K}( A ) := \ker h_\bc[A] = \ker ( c_0 {\bf 1}_{N \times N} + c_1
A + \cdots + c_{N - 1} A^{\circ ( N - 1 )} ).
\end{equation}

\noindent Then $\mathcal{K}(A) = \bigcap_{n \geq 0} \ker A^{\circ n}$,
and the extreme critical value is finite for all $A$:
\begin{equation*}
\mathcal{C}( h_\bc; z^M; A ) =
\sup_{\bu \in \mathcal{K}(A)^\perp \setminus \{ {\bf 0} \}} %
\frac{\bu^* A^{\circ M} \bu}%
{\bu^* \Bigl( \sum_{j = 0}^{N - 1} c_j A^{\circ j} \Bigr) \bu} \leq %
\mathcal{C}( \bc; z^M; N, \rho ), \qquad \forall A \in \bp_N(
\overline{D}( 0, \rho ) ) \setminus \{ 0 \}.
\end{equation*}

\noindent Moreover, the bound $\mathcal{C}( \bc; z^M; N, \rho )$ is
sharp, and is obtained as the supremum of the Rayleigh constant
$\mathcal{C}( h_\bc; z^M; A)$ as $A$ runs over the smaller set $\bp_N^1(
( \rho - \epsilon, \rho) )$ for any $\epsilon \in (0,\rho)$.
\end{proposition}

\noindent {\bf Sketch of proof.}
The first step is to show how the Schur polynomials $s_{\mu(M, N, j)}$ in
Theorem~\ref{Pjacobi-trudi} serve an additional purpose: they are
precisely the ``universal coefficients" involved in expressing $A^{\circ
M}$ as a combination of lower Hadamard powers, for any matrix $A$ and
over \textit{any} field $\F$. More precisely, if $A$ is an $N \times N$
matrix with entries in $\F$, and $\ba_1$, \ldots, $\ba_N$ are its
rows, then we first claim that
\begin{equation}\label{Emiracle2}
A^{\circ M} = \sum_{j = 0}^{N - 1} D_{M, j}( A ) A^{\circ j},
\end{equation}

\noindent where $D_{M, j}( A )$ is the diagonal matrix
$( {-1} )^{N - j - 1} \diag\bigl( s_{\mu( M, N, j )}( \ba_1 ), \ldots, 
s_{\mu( M, N, j )}( \ba_N ) \bigr)$. The claim follows by working with
distinct transcendental variables $s_1$, \ldots, $s_N$ and
solving the equation $V(\bu) {\bf s} = \bu^{\circ M}$ for ${\bf s}$,
where $V(\bu) := (u_i^{j-1})$ is the Vandermonde matrix, and ${\bf s} :=
(s_1, \ldots, s_N)^T$. The solution, via Cramer's rule, is given by:
$s_i = (-1)^{N - i} s_{\mu(M, N, i-1)}( \bu )$; now specialize to
$\bu = \ba_j^T$ for all $j$.

Having proved the claim, the second step is to show that $\mathcal{K}(A)
= \ker h_\bc[A] \subset \ker A^{\circ M}$ for all $M \geq 0$. This is
obvious if $0 \leq M < N$, while for $M \geq N$, we use~\eqref{Emiracle2}
to compute:
\[
h_\bc[A] \bv = 0\ \ \implies\ \ A^{\circ j} \bv = 0\ (0 \leq j < N)\ \
\implies\ \ \sum_{j=0}^{N-1} D_{M,j}(A) A^{\circ j} \bv = 0\ \
\implies\ \ A^{\circ M} \bv = 0.
\]

\noindent It follows that $\mathcal{K}(A) \subset \bigcap_{n \geq 0} \ker
A^{\circ n}$. The reverse inclusion is easy to show, as is (the equality
in) the next assertion. The subsequent inequality and the last sentence
in the result follow from Theorem~\ref{Tthreshold}.
\qed\medskip

It is also of interest to find a closed-form expression for
the generalized Rayleigh quotient $\mathcal{C}( h_\bc; z^M; A)$ for
a given matrix $A$. The following result provides two such expressions,
consequently revealing new and unexpected connections between
Rayleigh quotients and Schur polynomials.

\begin{proposition}
Fix integers $M \geq N \geq 1$ and positive scalars $c_0$, \ldots, $c_{N -
1} > 0$. Then,
\begin{equation}\label{Especrad}
\mathcal{C}( h_\bc; z^M; A ) =
\varrho( h_\bc[A]^{\dagger/2} A^{\circ M} h_\bc[A]^{\dagger/2}),
\qquad \forall A \in \bp_N(\C) \setminus \{ 0 \},
\end{equation}

\noindent where $C^{\dagger/2}, \varrho(C)$ denote the principal square
root of the Moore-Penrose inverse of $C$, and the spectral radius of $C$,
respectively. For instance, if $A = \bu \bu^*$
with $\bu$ having distinct coordinates, then
\begin{equation}\label{Erankone}
\mathcal{C}( h_\bc; z^M; \bu \bu^* )
= (\bu^{\circ M})^* h_\bc[ \bu \bu^*]^\dagger \bu^{\circ M}
= \sum_{j=0}^{N-1} \frac{|s_{\mu(M,N,j)}(\bu)|^2}{c_j}.
\end{equation}
\end{proposition}

\noindent \textbf{Sketch of proof.}
The proof of~\eqref{Especrad} uses the theory of Kronecker normal forms
and Rayleigh quotients, and is omitted for brevity. To show the first
equality in~\eqref{Erankone}, set $\bv := h_\bc[\bu \bu^*]^{\dagger/2}
\bu^{\circ M}$. Then standard computations show that
$\mathcal{C}(h_\bc; z^M; \bu \bu^*) = \varrho(\bv \bv^*) =
\bv^* \bv = (\bu^{\circ M})^* h_\bc[\bu \bu^*]^\dagger \bu^{\circ M}$.

Finally, we show that the last equality in~\eqref{Erankone} holds more
generally, for any rank-one matrix $A = \bu \bv^T$, where $\bu, \bv$ are
vectors with distinct coordinates in any field $\F$. Indeed,
notice by the proof of~\eqref{Ecauchy-binet1} that
\[
h_\bc[ \bu \bv^T ] = %
X( \bu, \bn_{\min} ) \diag( c_0, \ldots, c_{N - 1} ) %
X( \bv, \bn_{\min})^T, \ \ \text{where} \ \ %
\bn_{\min} := ( 0, 1, \ldots, N - 1 ).
\]
Moreover, $X( \bu, \bn_{\min} )$ is precisely $V(\bu)$, the Vandermonde
matrix $(u_i^{j - 1})$. Now the proof of Proposition~\ref{Crayleigh}
shows that
$V(\bu)^{-1} \bu^{\circ M} = ((-1)^{N - j - 1} %
s_{\mu(M,N,j)}(\bu))_{j=0}^{N-1} =: {\bf s}_N(\bu)$,
say.  Hence,
\begin{align*}
(\bv^{\circ M})^T h_\bc[\bu \bv^T]^{-1} \bu^{\circ M}
= &\ (V(\bv)^{-1} \bv^{\circ M})^T \diag(c_0, \ldots, c_{N-1})^{-1}
(V(\bu)^{-1} \bu^{\circ M})\\
= &\ {\bf s}_N(\bv)^T \diag(c_0, \ldots, c_{N-1})^{-1} {\bf s}_N(\bu)
= \sum_{j=0}^{N-1} \frac{s_{\mu(M,N,j)}(\bu) s_{\mu(M,N,j)}(\bv)}{c_j}.
\hspace*{2mm}\square
\end{align*}

Equation~\eqref{Erankone} provides an alternate explanation of how and
why Schur polynomials occur in the extreme critical value $\mathcal{C}(
\bc; z^M; N, \rho )$, by considering the matrices in
$\bp_N^1((0,\rho))$. Having carried out step (I) in trying to prove
Theorem~\ref{Tthreshold} by an alternate approach (see the previous
page), a natural question is to ask if it is possible to maximize the
function $\Psi_{\bc, M} : A \mapsto \mathcal{C}( h_\bc, z^M; A)$ to
obtain $\mathcal{C}( \bc; z^M; N, \rho)$, as in step (II). Observe
by~\eqref{Especrad} that the spectral map $\Psi_{\bc, M}$ is continuous
on the open dense subset of the cone given by $\det h_\bc[A] \neq 0$.
However, $\Psi_{\bc, M}$ turns out to \textit{not} be continuous
on all of $\bp_N( \overline{D}( 0, \rho ) )$, or even on
$\bp_N^1([0,\rho])$. Specifically, it is not continuous at the matrix $A
= \rho \one{N}$. This spectral discontinuity phenomenon warrants further
exploration, and is to be the subject of future
work~\cite{BGKP-fixeddim2}.

\section{Stratification of the cone, and the simultaneous kernels}

In the final section, we take a closer look at the simultaneous kernel
$\mathcal{K}(A)$ defined in~\eqref{Esimult}. As we now discuss, this
space crucially depends on a canonical block decomposition of the matrix
$A$. We begin by isolating this refined structure. Consider the
following two examples:
\[
A_1 := \begin{pmatrix} 5 \cdot {\bf 1}_{a \times a} & B\\ B^* & 2
\cdot {\bf 1}_{b \times b} \end{pmatrix} \in \bp_{a+b}(\C), \qquad
A_2 := \begin{pmatrix} 5 & -5 & u_1\\ -5 & 5 & u_2\\
\overline{u_1} & \overline{u_2} & 2 \end{pmatrix} \in \bp_3(\C).
\]

\noindent Given the positivity of $A_1, A_2$, one can show that all
entries of $B$ are equal, while $u_1 = - u_2$. In fact, if
entries in each diagonal block of a positive semidefinite matrix lie in a
$G$-orbit for some subgroup $G \subset \C^\times$, this imposes
constraints on the off-diagonal blocks. This is distilled into the
following result.

\begin{theorem}
Fix a subgroup $G \subset \C^\times$, an integer $N \geq 1$,
and a non-zero matrix $A \in \bp_N( \C )$.
There exists a partition $\pi^G(A) := \{ I_1, \ldots, I_k \}$ of
$\{ 1, \ldots, N \}$ (unique up to relabelling), satisfying:
\begin{enumerate}
\item Each diagonal block $A_{I_j}$ of $A$ is a submatrix with rank at
most one.

\item The entries of each diagonal block $A_{I_j}$ lie in a single
$G$-orbit.

\item The diagonal blocks $A_{I_j}$ of $A$ satisfying (1), (2) have
maximal size.
\end{enumerate}
In this case, each off-diagonal block of $A$ also has rank at most one,
with all its entries in a single $G$-orbit.
\end{theorem}

\noindent For instance, for two choices of the group $G$ the partition
$\pi^G(A)$ is easily interpreted:

1. $G = \{ 1 \}$, in which case all entries in a diagonal (or
off-diagonal) block of $A$ are equal.

2. $G = S^1$, in which case all entries in a diagonal (or off-diagonal)
block of $A$ are equal in modulus.\medskip

\noindent {\bf Sketch of proof.}
Suppose $\{ I_1, \ldots, I_k \}$ is any partition of $\{ 1, \ldots, N \}$
satisfying conditions (1), (2). Let $1 \leq i \neq j \leq k$, and $1 \leq
l < l' < m \leq N$, with $l$, $l' \in I_i$ and $m \in I_j$. Consider the
submatrix
$B := A_{\{ l, l', m \}} = \begin{pmatrix}
a & a g & b \\
a \overline{g} & a | g |^2 & c \\
\overline{b} & \overline{c} & d
\end{pmatrix}$,
where $a, d \geq 0,\ g, \overline{g} \in G,\ b, c \in \C$.
We claim that $c \in b \cdot G$, and that the minor
$\begin{pmatrix} a & b \\ a \overline{g} & c \end{pmatrix}$
is singular. This is because
$0 \leq \det B = %
-a ( | c |^2 + | b |^2 | g |^2 - 2 \Re( \overline{b} c g ) ) = %
-a | c - b \overline{g} |^2$.
Hence either $a = 0$, in which case $b = c = 0$, by the positivity
of~$B$, or $c = b \overline{g}$. The proof repeatedly uses computations
along similar lines, to show that there exists 
$C \in \bp_k(\C)$ with ${\rm rank} (C) = {\rm rank} (A)$,
and vectors $\bu_j \in \C^{|I_j|}$ with entries in a single $G$-orbit,
such that $A_{I_i \times I_j} = c_{ij} \bu_i \bu_j^*$, for all $1 \leq
i,j \leq k$.
\qed\medskip

Denote by $(\Pi_N, \prec)$ the poset of all partitions of $\{ 1, \ldots, N
\}$ under refinement. Then one has the partition map $\pi^G : \bp_N(\C)
\to \Pi_N$, sending $0$ to $\{ \{ 1, \ldots, N \} \}$ and all other
matrices $A$ to $\pi^G(A)$. Define $\cals^G_\pi$ to be the fiber of this
map:
\begin{equation}
\cals^G_\pi := \{ A \in \bp_N(\C) : \pi^G(A) = \pi \},
\qquad \forall \pi \in \Pi_N.
\end{equation}

\begin{cor}
Fix a subgroup $G \subset \C^\times$. The sets $\cals^G_\pi$
form a Schubert cell-type stratification of the cone:
\begin{equation}
\bp_N(\C) = \bigsqcup_{\pi \in \Pi_N} \cals^G_\pi, \qquad
\overline{\cals^G_\pi} = \bigsqcup_{\pi' \prec \pi} \cals^G_{\pi'}, \qquad
\forall N \geq 1,\ \pi \in \Pi_N.
\end{equation}

\noindent Moreover, every $A \in \bp_N(\C)$ has rank at most
$|\pi^{\C^\times}(A)|$.
\end{cor}

The stratification of the cone $\bp_N( \C )$ is noteworthy in that the
generalized Rayleigh quotient map $\Psi_{\bc,M}$ (defined in
Section~\ref{Srayleigh}) is discontinuous at the point $\rho \one{N}$ as
one is jumping across strata $\cals^{\{ 1 \}}_\pi$.

Motivated by Proposition~\ref{Crayleigh}, a precise description of the
simultaneous kernel $\mathcal{K}(A) = \bigcap_{n \geq 0} \ker A^{\circ
n}$ is in order. It turns out that the map $A \mapsto \mathcal{K}(A)$
depends crucially (and solely) on the stratification.

\begin{theorem}\label{Tsimult}
The simultaneous kernel map $A \mapsto \mathcal{K}(A)$ is constant on
each stratum $\mathcal{S}^{\{ 1 \}}_\pi$, i.e.,
\[
\mathcal{K} : \bp_N(\C) \longrightarrow \Pi_N \longrightarrow 
\bigsqcup_{r = 0}^{N - 1} \Gr( r, \C^N )
\]

\noindent sends every matrix $A \in \mathcal{S}^{\{ 1 \}}_\pi$ to a fixed
subspace
\[
\mathcal{K}_\pi := \ker \oplus_j {\bf 1}_{I_j \times I_j} =
\oplus_j \ker {\bf 1}_{I_j \times I_j} \in \Gr( N - | \pi |, \C^N ),
\]

\noindent where $\pi = \{ I_j \}$, and $\Gr( r, \C^N )$ is the complex
Grassmann manifold of $r$-dimensional subspaces of~$\C^N$.
\end{theorem}

\noindent The proof of this result is fairly involved, and we refer the
reader to \cite[Section 5]{BGKP-fixeddim} for details.

We conclude with the following surprising consequence of
Theorem~\ref{Tsimult}: as $A$ runs over the uncountable set of matrices
in $\bp_N(\C)$, the set of simultaneous kernels $\mathcal{K}(A) =
\bigcap_{n \geq 0} \ker A^{\circ n}$ is, nevertheless, a
finite set of subspaces of $\C^N$, indexed by $\Pi_N$.
This is in stark contrast to the situation for the usual matrix
powers, in which case $\bigcap_{n \geq 1} \ker A^n = \ker A$ can vary
over an uncountable set of subspaces of $\C^N$.

Other ramifications of this work, as well as complete proofs can be found
elsewhere \cite{BGKP-fixeddim, BGKP-fixeddim2}.


%
%


\end{document}